\documentclass[12pt]{article}

\usepackage{theorem,amssymb,epsfig}

\advance \topmargin by -\headheight
\advance \topmargin by -\headsep
\evensidemargin \oddsidemargin
\author{J.-P. Allouche \\
CNRS, Math., \'Equipe Combinatoire et Optimisation \\
Universit\'e Pierre et Marie Curie, Case 189 \\
4 Place Jussieu \\
F-75252 Paris Cedex 05 \\
France \\
{\tt allouche@math.jussieu.fr}
\and
M. Mend\`es France \\
Universit\'e Bordeaux I \\
Math\'ematiques \\
F-33405 Talence Cedex \\
France \\
{\tt michel.mendes-france@math.u-bordeaux1.fr} \\
}

\title{Hadamard grade of power series}

\date{ }

\def \proof{\bigbreak\noindent{\it Proof.\ \ }}

\def \endpf{{\ \ $\Box$ \medbreak}}

\newtheorem{theorem}{Theorem}

\newtheorem{corollary}{Corollary}

\theorembodyfont{\rm}

\newtheorem{remark}{Remark}
\newtheorem{example}{Example}
\newtheorem{definition}{Definition}

\begin{document}

\maketitle

\begin{center}
{\it \`A Vera S\'os pour son $80^e$ anniversaire, avec admiration et amiti\'e}
\end{center}

\medskip

\begin{abstract}
The Hadamard product of two power series $\sum a_n z^n$ and 
$\sum b_n z^n$ is the power series $\sum a_n b_n z^n$. We define 
the (Hadamard) grade of a power series $A$ to be the least number 
(finite or infinite) of algebraic power series, the Hadamard product 
of which equals $A$. We study and discuss this notion.

\medskip
{\scriptsize


}

\end{abstract}

\section{Introduction}

The Hadamard product of two power series is simply the power series
obtained by termwise multiplication of their coefficients. Since
the seminal paper of Hadamard \cite{Hadamard} this product occurred
in various, sometimes unexpected, fields (look for example at the
introduction of \cite{Bragg}). We revisit the question of what can be
said about the Hadamard product of ``simple'' series. We introduce
in particular the closure under finite Hadamard products of the set
of algebraic power series.

\section{The early work of Hadamard}

\begin{definition}
Let $A(z) := \sum_{n \geq 0} a_n z^n$ and $B(z) := \sum_{n \geq 0} b_n z^n$ 
be two power series with coefficients in some commutative field $K$. 
Their {\em Hadamard product} (sometimes also called {\em child product}, see 
\cite[Section~8, p.\ 251]{vdP}, {\em Schur product}, or {\em quasi-inner product},
see \cite[p.\ 36]{Bragg}) is the power series  $A * B (z)$ defined by
$$
A * B (z) := \sum_{n \geq 0} a_n b_n z^n.
$$
\end{definition}

Note that this definition makes sense either by considering series
$A$ and $B$ as formal power series, or, if $K$ is a subfield of the 
field of complex numbers ${\mathbb C}$, as series possibly converging 
in a neighborhood of $0$ (if $A$ and $B$ converge in a neighborhood of
$0$, then obviously $A*B$ also converges in some neighborhood of $0$).
If $A$ and $B$ are two functions, we let $A * B$ denote the Hadamard
product of their respective power series, provided they exist.

\subsection{Singularities of Hadamard products}

Starting from the easy observation that for $a$ and $b$ nonzero
$$
\frac{1}{a-z} * \frac{1}{b-z} = \frac{1}{ab-z}
$$
one may ask whether the singularities of the Hadamard product of
two series with coefficients in ${\mathbb C}$ are related to the
product of singularities of these series. Hadamard realized that 
this is indeed the case.

\begin{theorem}[Hadamard \cite{Hadamard}]\label{sing-hadamard}
Let $A(z) = \sum a_n z^n$ and $B(z) = \sum b_n z^n$ be two power
series holomorphic in a neighborhood of $0$. Let $\sigma_A$
and $\sigma_B$ be their respective sets of singularities. Then
the singularities of $A * B$ are among the values $\alpha \beta$,
where $\alpha \in \sigma_A$ and $\beta \in \sigma_B$. In other words
$\sigma_{A*B} \subset \sigma_A \cdot \sigma_B := 
\{\alpha \beta, \ \alpha \in \sigma_A, \ \beta \in \sigma_B\}$.
\end{theorem}

The proof, which is by no means easy, involves the study of the integral
$$
\oint A(u)B\left(\frac{z}{u}\right) \frac{{\rm d}u}{u} =\sum a_n b_n z^n
$$
extended to a circle around the origin. The holomorphy of $A(z)$ and $B(z)$
in a neighborhood of $0$ is therefore an essential hypothesis for the result to
hold. And indeed in Section~\ref{conv} below we show what may happen if say
$\limsup_{n +\infty} |a_n|^{1/n} = \infty$.

\begin{remark}
Other papers on singularities of Hadamard product may be of interest for
the reader: see in particular \cite{Borel, Faber, FFK}.
\end{remark}

\subsection{A counter-example: the Euler series}\label{conv}

The hypothesis that both $A$ and $B$ converge in some neighborhood of the 
origin is necessary for the validity of Theorem~\ref{sing-hadamard}. 
Indeed consider the Euler $E(z)$ formal power series defined by
$$
E(z) := \sum_{n \geq 0} (-1)^n n! z^n.
$$
This series, considered as a complex series, diverges everywhere except 
at $0$. As will be explained below, $E(z)$ has the integral representation 
$I(z)$ where
$$
(\#) \hskip 2truecm I(z) = \int_0^{\infty} \frac{e^{-u}}{1+zu} \, {\rm d}u =
- \frac{1}{z} \, e^{1/z} \log \frac{1}{z} + S\left(\frac{1}{z}\right)
$$
with 
$$
S(y) := -y e^y \left(\gamma + \sum_{n \geq 1} \frac{(-1)^n y^n}{n. n!}\right)
$$
and $\gamma$ is the Euler constant.
The integral above is a function well-defined in the complex plane
except on the real negative axis $\{z = s+it, \ t = 0, \ s < 0\}$.

\medskip

\noindent
Now let $B(z) := e^z$. 
We have $E * B (z) = \sum_{n \geq 0} (-1)^n z^n =\frac{1}{1+z}$ which has
only one singularity at $z = -1$, while $e^z$ does not have any singularity.

\medskip

In order to demystify the integral representation $E(z) = I(z)$, compute 
the derivatives of $I(z)$:
$$
I^{(n)}(z) = (-1)^n \int_0^{\infty} \frac{n!e^{-u}}{(1+zu)^{n+1}} \, {\rm d}u
$$
which reduces at $z=0$ to $I^{(n)}(0) = (-1)^n (n!)^2$.
Therefore the formal Taylor series for $I(z)$ is 
$$
\sum_{n \geq 0} (-1)^n n! z^n
$$
which indeed coincides with $E(z)$.

\medskip

The last part of identity $(\#)$ is an illuminating computation (see Hardy's 
book \cite[p.~26--27]{Hardy}) that describes precisely the singularities 
of $E(z)$.

\begin{theorem}[\cite{Hardy}]
The integral rerpresentation of the Euler series 
$\sum_{n \geq 0} (-1)^n n! z^n$ can be continued for all $z$
to a many-valued function with an infinity of branches differing by integral
multiples of $\, 2i\pi z^{-1}e^{1/z}$. One branch tends to $1$ 
as $z$ tends to $0$ through positive values. More precisely
$$
I(z) = - \frac{1}{z} \, e^{1/z} \log \frac{1}{z} + S\left(\frac{1}{z}\right)
$$
where 
$$
S(y) := -y e^y \left(\gamma + \sum_{n \geq 1} \frac{(-1)^n y^n}{n. n!}\right)
$$
and $\gamma$ is the Euler constant. $I(z)$ lives on the same Riemann surface as 
$\log z$. In particular, for $z=1$, {\em [Euler and Hardy write]}
$$
\sum_{n \geq 0} (-1)^n n! 
= - e \left(\gamma + \sum_{n \geq 1} \frac{(-1)^n}{n n!}\right) = .5963...
$$
\end{theorem}

\proof (This proof is attributed to Euler by Hardy in \cite{Hardy}).

Suppose $z = x$ is a positive real number. Performing the change of
variable $t = \frac{x}{1+xu}$, we have
$$
I(x) = \int_0^{\infty} \frac{e^{-u}}{1+xu} \, {\rm d}u 
= - \frac{1}{x} \, e^{1/x} \, {\rm li}(e^{-1/x})
$$
where ${\rm li}(\xi) := \displaystyle\int_0^{\xi} \frac{{\rm d}s}{\log s}$, 
for $\xi \in [0, 1)$, is the {\em logarithmic integral} ({\em logarithm-integral\,}
in Hardy's terminology). Now
$$
\begin{array}{lll}
- {\rm li}(e^{-y}) &=& \displaystyle\int_y^{\infty} \frac{e^{-u}}{u} \, {\rm d}u
= \int_1^{\infty} \frac{e^{-u}}{u} \, {\rm d}u 
  - \int_0^1 \frac{1-e^{-u}}{u} \, {\rm d}u - \int_1^y \frac{{\rm d}u}{u}
  + \int_0^y \frac{1 - e^{-u}}{u} \, {\rm d}u \\
&=& - \gamma - \log y - \displaystyle\sum_{n \geq 1} \frac{(-1)^n y^n}{n. n!}
\end{array}
$$
which concludes the proof. \endpf

\section{Digression: elementary optics}\label{digression}

Consider $N$ transparent plates one against the other with
respective thicknesses $a_1, a_2, \ldots, a_N$ ($a_n > 0$)
and indices of refraction $n_1, n_2, \ldots, n_N$ ($n_k > 1$).

\centerline{\epsfig{figure=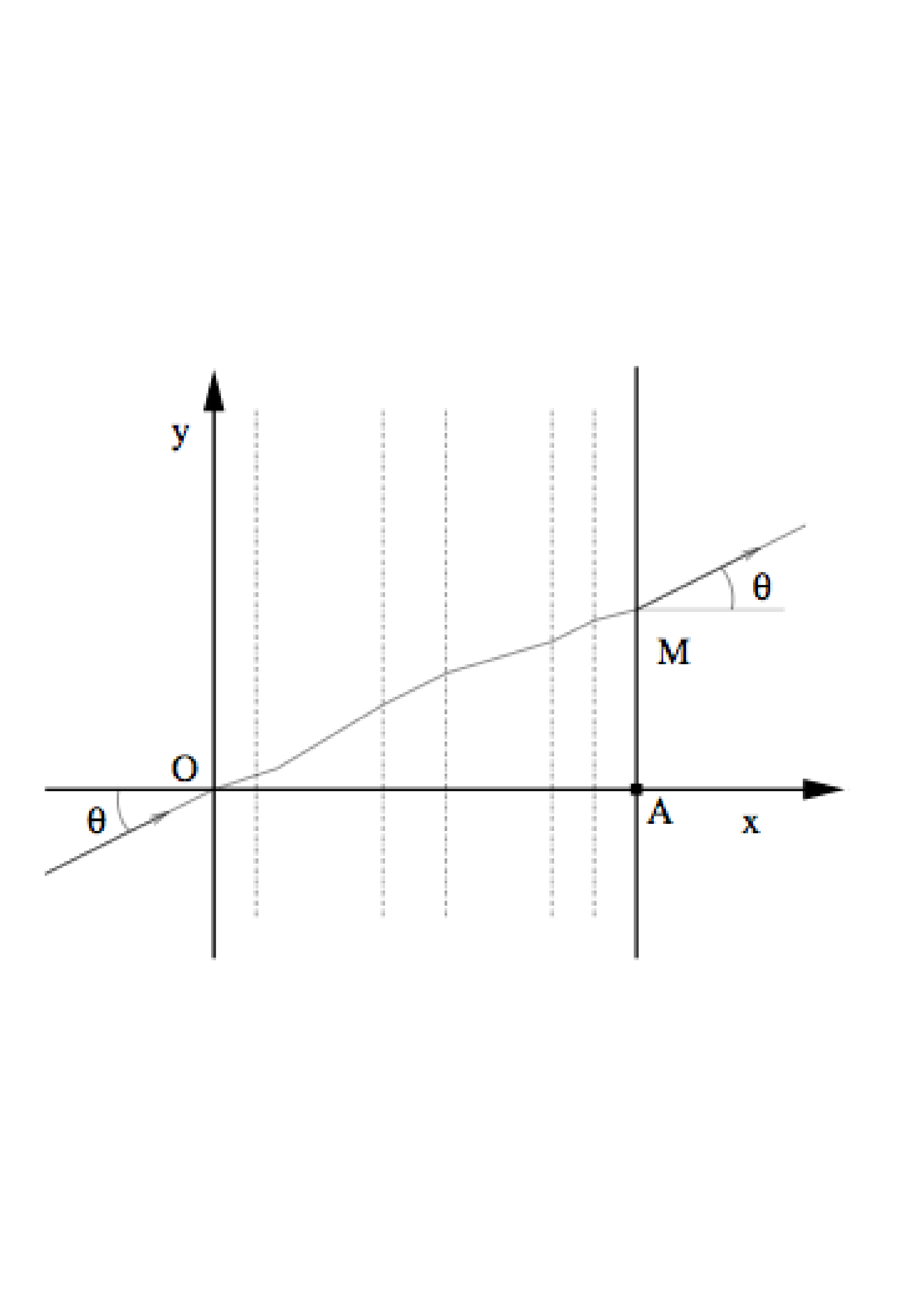,height=9cm}}

An incoming ray hits the system at the origin with an incident
angle $\theta$ and reappears after $N$ refractions at $M$ on
the opposite side. Let $\theta_k$ be the angle at which the luminous
path meets the $k$th plate. The Snell-Descartes law states that
$n_k \sin \theta_k$ is independent of $k$. An easy computation shows
that the length $AM = \varphi(\theta)$ is given by
$$
\varphi(\theta) = 
\sum_{1 \leq k \leq N} \frac{a_k \sin \theta}{\sqrt{n_k^2 - \sin^2 \theta}} \cdot
$$
Put $z := \sin \theta$ and $\psi(z) := \varphi(\theta)$, thus
$$
\psi(z) = \sum_{1 \leq k \leq N} \frac{a_k z}{\sqrt{n_k^2 - z^2}} \cdot
$$
In \cite{MFS} the second named author and A. Sebbar notice that
$$
\psi(z) = 
\frac{z}{\sqrt{1 - z^2}} \ * \sum_{1 \leq k \leq N} \frac{a_k n_k z}{n_k^2 - z^2} \cdot
$$
If none of the $a_k$ vanish and if the $n_k$ are pairwise distinct, then $\psi$
is an algebraic function of degree $2^n$. The above Hadamard representation is
somewhat simpler than the original $\psi$ since one factor is rational whereas 
the other is quadratic.

The above identity is even valid for $N = \infty$, $a_k$ and $n_k \in {\mathbb C}$,
provided the series converges. The second factor can often be summed, especially
if $a_k$ and $n_k$ are rational functions of $k$. Choose for example $n_k = (2k+1)$
and $a_k = (2k+1)^{-1}$. Then the identity reads:
$$
\sum_{k \geq 0} \frac{z}{(2k+1)\sqrt{(2k+1)^2 - z^2}} =
\frac{\pi}{4} \frac{z}{\sqrt{1 - z^2}} * \tan \frac{\pi z}{2}\cdot
$$
This could be thought as of a quasi-summation of the lefthand side series.

Let us test the formula. For $z \to 0$ it reduces to
$$
\sum_{k \geq 0 } \frac{1}{(2k+1)^2} = \frac{\pi^2}{8}\cdot
$$
More generally, put $\tan X := \displaystyle\sum_{j \geq 0} t_{2j+1} X^{2j+1}$.
Observe that
$$
\frac{1}{(2k+1)\sqrt{(2k+1)^2 - z^2}} = 
\frac{1}{(2k+1)^2} 
\sum_{j \geq 0} (-1)^j {-\frac{1}{2} \choose j} \frac{z^{2j}}{(2k+1)^2j}\cdot
$$
Then the above Hadamard product leads to the well known equality
$$
\sum_{k \geq 0} \frac{1}{(2k+1)^{2j+2}} = \frac{\pi^{2j+2}}{2^{2j+3}} t_{2j+1}.
$$

The example $a_k = (2k+1)^{-1}$ and $n_k = (2k+1)$ corresponds to an optical system
with infinitely many thinner and thinner plates but with no physical meaning since
$n_k = (2k+1)$ tends to infinity. In real materials the highest indices seem not to go
beyond $4$ or $5$.

\bigskip

One last remark before closing this section. The identity showing that
$\psi$ is a Hadamard product of two simple functions can easily be extended.
Let $H(z) := \displaystyle\sum_{j \geq 0} h_j z^{2j+1}$ be any odd function holomorphic
around the origin ($|z| < R$). Then, for all complex nonzero $a_k$ and $n_k$,
and for $N \leq +\infty$,
$$
\sum_{1 \leq k \leq N} a_k H\left(\frac{z}{n_k}\right) =
H(z) \ * \sum_{1 \leq k \leq N} \frac{a_k n_k z}{n_k^2 - z^2}
$$
provided $|z| < R \, \displaystyle\inf_k |n_k|$. The proof consists in expanding both 
sides of the identity into Taylor series, just as in the case $H(z) = z (1 - z^2)^{-1/2}$.
Needless to say that similar identities hold if $H$ is even, or neither odd nor even.

\section{Recalling some definitions and results}

In what follows we will need the definition of {\em $D$-finite\,} power series
and the definition of {\it automatic\,} sequences.

\begin{definition}[\cite{Stanley}]
A power series is {\em $D$-finite} (or {\em differentially finite} or {\em holonomic}) 
if it satisfies a linear differential equation with polynomial coefficients.
\end{definition} 

\begin{example}
The power series $F := \sum_{n \geq 0} \frac{z^n}{n!}$ is $D$-finite, since it satisfies
$F' - F =0$. The series $\sum_{n \geq 0} \sqrt{n} z^n$ is not $D$-finite
\cite{Gerhold}. The series $\sum_{n \geq 1} \log n \ z^n$ and $\sum_{n \geq 1} p_n z^n$
where $p$ is the $n$-th prime are not $D$-finite \cite{FGS}.
A nice example of a non-$D$-finite power series is $\sum_{n \geq 1} \zeta(2n+1) z^n$ 
where $\zeta$ is the Riemann function (see \cite{BGKL}).

\end{example}

\begin{definition}
Let $q$ be an integer $\geq 2$. A sequence $(a_n)_{n \geq 0}$ is {\em $q$-automatic} 
if the $q$-kernel of $(a_n)_{n \geq 0}$, i.e., the set of subsequences
$\{(a_{q^k n + j})_{n \geq 0}, \ k \geq 0, \ j \in [0, q^k - 1]\}$, is finite.
\end{definition}

\section{The Hadamard grade of a power series}

\begin{definition}
A power series $A(z) := \sum a_n z^n$ with coefficients in a commutative field $K$ is
said to have {\em finite Hadamard grade over $K$} if there exist finitely many power 
series $B_1(z)$, $B_2(z)$, ..., $B_k(z)$ that are algebraic (over $K(z)$) such that 
$A = B_1 * B_2 * \cdots * B_k$. The least $k$ having this property is called the 
{\em Hadamard grade over $K$} of $A$.
\end{definition}

We give some first properties of power series with finite Hadamard grade.

\begin{theorem}

\ { }

\begin{itemize}

\item The Hadamard product of two power series having finite Hadamard grade has finite 
Hadamard grade.

\item If $A$ is a power series that is algebraic over $K(z)$, then $A$ has Hadamard 
grade $1$.

\item If $A$ is a power series with finite Hadamard grade, then $A$ is {\em $D$-finite}.

\end{itemize}

\end{theorem}

\proof The first two assertions are clear. For the third, recall that a power series is 
{\em $D$-finite} (or {\it differentially finite} or {\em holonomic}) if it satisfies a 
linear differential equation with polynomial coefficients (see \cite{Stanley}); 
also recall that any algebraic power series is $D$-finite and that the Hadamard product 
of two $D$-finite power series is also $D$-finite (see \cite{Stanley}).  \endpf

\begin{remark}
The $D$-finiteness of the Hadamard product of two $D$-finite power series can be found
in the paper \cite{Stanley}, but as noted by Stanley this is already in the paper 
of Jungen \cite{Jungen} who indicates that this result can be found in the unpublished 
manuscripts of Hurwitz (available at \'Ecole Polytechnique de Zurich). Jungen also gives 
the same characterization as Stanley's for the coefficients of $D$-finite power series.
\end{remark}

The next result shows that in positive characteristic the definition of the grade is 
not really pertinent (that is why in the next sections we will focus on the case where 
$K$ is the field of complex numbers ${\mathbb C}$ or the field of rational numbers
${\mathbb Q}$).

\begin{theorem}\label{charpos}
A power series on a field of positive characteristic has finite Hadamard grade if and
only if it is algebraic.
\end{theorem}

\proof Only one direction need to be proven. This was done by Furstenberg 
\cite{Furstenberg} for a finite field. Fliess \cite{Fliess} noted that the proof 
still works for a perfect field of positive characteristic. Deligne \cite{Deligne} 
gave a general proof. \endpf

\begin{remark}
For more around algebraic power series the reader can look for example at 
\cite{Allouche89, Griffiths1, Griffiths2, vdP} and the references therein. 
\end{remark}

\section{Hadamard grade of formal power series with complex or with rational coefficients}

\subsection{The base field is ${\mathbb C}$}

    From now on we will restrict ourselves to the case where the field $K$ is either 
${\mathbb C}$ or ${\mathbb Q}$. We begin with a theorem stating  that a minimal 
Hadamard decomposition of an irrational power series into finitely many algebraic power 
series only involves algebraic irrational power series.

\begin{theorem}
Let $A(z)$ be a power series with complex coefficients. If $A$ is irrational
and has finite Hadamard grade $d$, then in any Hadamard decomposition
$A = B_1 * B_2 * \ldots * B_d$ with algebraic $B_i$'s, all the $B_i$'s are irrational.
\end{theorem}

\proof It suffices to use the minimality of $d$ and a result of Jungen \cite{Jungen} 
stating that the Hadamard product of a rational power series and an algebraic power 
series is algebraic.
\endpf

\medskip

The following result is classical. 

\begin{theorem}
Let $F =\sum_{n \geq 0} a_n z^n$ be a formal power series with coefficients in
the complex field $\,{\mathbb C}$ that is algebraic over $\,{\mathbb C}(z)$. Then 
the radius of convergence of $F$ is $>0$ and $F$ has finitely many singularities 
in its disk of convergence.
\end{theorem}

Using this theorem and Hadamard's Theorem \ref{sing-hadamard} yields the following 
corollary.

\begin{corollary}\label{radius}
Let $A(z)$ be a formal power series with coefficients in ${\mathbb C}$ which has
finite Hadamard grade. Then $A$ has a positive radius of convergence. Furthermore
$A$ has only finitely many singularities in its disk of convergence.
\end{corollary}

There is a simple necessary condition for a complex power series to be of finite 
Hadamard grade in terms of diagonals of rational functions.

\begin{theorem}
If $A$ is a complex power series with finite Hadamard grade $d$, then $A$ is the
(complete) diagonal of a complex rational power series with $2d$ variables.
\end{theorem}

\proof The proof is done by induction on $d$. The case $d=1$ can be found in the
paper of Furstenberg \cite{Furstenberg}. Now if $A_1, A_2, \ldots, A_d$ are complex 
algebraic power series, then there exist a two-variable rational power series
$B(X,Y) := \sum b_{k,\ell} X^k Y^{\ell}$ and a rational power series
$C(Z_1, Z_2, \ldots, Z_{2d-2}) := \sum c_{j_1, j_2, \ldots, j_{2d-2}} Z_1^{j_1} Z_2^{j_2} 
\ldots Z_{2d-2}^{j_{2d-2}}$ such that $A_d(z) = \sum b_{k, k} z^k$ and
$A_1(z) * A_2(z) * \ldots A_{d-1}(z) = \sum c_{j, j, \ldots, j} z^j$.
The (ordinary Cauchy) product of the series $B$ and $C$, i.e.,
$B(X,Y) C(Z_1, Z_2, \ldots, Z_{2d-2})  =
\sum b_{k,\ell} c_{j_1, j_2, \ldots, j_{2d-2}} X^k Y^{\ell} Z_1^{j_1} Z_2^{j_2} 
\ldots Z_{2d-2}^{j_{2d-2}}$ is a rational power series whose complete diagonal
$\sum b_{j,j} c_{j, j, \ldots, j} z^j$ is equal to $A_1(z) * A_2(z) * \ldots A_d(z)$.
\endpf

\subsection{The base field is ${\mathbb Q}$}

What can be said of a power series with rational coefficients and finite Hadamard
grade over ${\mathbb Q}$? We begin with a simple remark.

\begin{remark}
The family of power series with rational coefficients and with finite Hadamard grade 
over ${\mathbb Q}$ is countable.
\end{remark}

Our next result links finite Hadamard grade over ${\mathbb Q}$ and automaticity
for the case of power series with rational coefficients (for more about
{\em automatic sequences} the reader can look at \cite{AS}).

\begin{theorem}
Let $A(z) := \sum a_n z^n$ be a power series with rational coefficients.
If $A$ has finite Hadamard grade over the field of rational numbers ${\mathbb Q}$,
then there are only finitely many primes dividing at least one of the denominators
of the $a_n$'s. Furthermore for all primes $p$ outside these finitely many primes, the 
series $\overline{A}(z) := \sum (a_n \bmod p) z^n$ is algebraic over ${\mathbb F}_p(z)$.
In particular the sequence $(a_n \bmod p)_{n \geq 0}$ is $p$-automatic. Even more,
for these primes $p$ and for any $r \geq 1$, the sequence $(a_n \bmod p^r)_{n \geq 0}$ 
is $p$-automatic.
\end{theorem}

\proof Write $A = B_1 * B_2 * \cdots * B_d$ where the $B_i$'s are algebraic over 
${\mathbb C}(z)$. From Eisenstein's theorem \cite{Eisenstein} (see Remark~\ref{Eis} 
below) the denominators of the coefficients of each $B_i$ have only finitely many 
distinct prime factors. Let ${\cal P}$ be the (finite) set of primes
that divide the denominator of at least one $B_i$. For each prime $p \notin {\cal P}$
one can write $(A \bmod p) = (B_1 \bmod p) * (B_2 \bmod p) * \cdots * (B_d \bmod p)$.
Each $(B_i \bmod p)$ is algebraic over ${\mathbb F}_p(z)$, so is their Hadamard
product (Theorem~\ref{charpos}). Actually for any $r \geq 1$ the sequence
$(a_n \bmod p^r)_{n \geq 0}$ is $p$-automatic: this was proved by Christol
(see \cite{Christol1, Christol2}, see also \cite{DL}). \endpf

\begin{remark}\label{Eis}

\ { }

\begin{itemize}

\item The Eisenstein theorem states: {\it Let $\sum a_n z^n$ be an algebraic
power series with rational coefficients. Then there exist two integers 
$A$ and $C$ such that, for all $n \geq 0$, the number $C A^n a_n$ is an integer}.
This theorem was actually proved by Heine \cite{Heine}.

\item Christol (see, e.g., \cite{Christol2}) calls {\it globally bounded} any series
$f = \sum a_n z^n$ with rational coefficients such that $f$ has a positive radius of
convergence (when seen as a power series on ${\mathbb C}$) and such that there exist
$\alpha, \beta \in {\mathbb Q}$ such that $\alpha f(\beta z)$ belongs to ${\mathbb Z}[[z]]$.
He shows that that any algebraic power series $\sum a_n z^n$ is {\em globally automatic}, 
i.e., for all but finitely many primes $p$ and for all integers $r \geq 2$, the sequence 
$(a_n \bmod p^r)_{n \geq 0}$ is $p$-automatic.

\end{itemize}

\end{remark}

\section{Examples, Counterexamples, and Questions}

Using the results of the previous sections we can give examples and counterexamples
of power series having finite Hadamard grade: we will restrict to power series {\bf with
rational coefficients}, and consider their {\bf Hadamard grade over ${\mathbb Q}$},
although some of the statements and questions below could be formulated for the Hadamard
grade over $K$ of power series with coefficients in $K$, where $K$ is any field of 
characteristic zero.

\begin{itemize}

\item[(a)] {\it The power series $e^z$, $\log(1+z)$, $\sum \pm z^n$ where the sequence
$\pm 1$ is not ultimately periodic, all have infinite Hadamard grade.}

\proof The denominators of the power series $e^z$ and $\log(1+z)$ are divisible
by infinitely many primes. If the sequence $\pm$ is not ultimately periodic,
then the series $\sum \pm z^n$ admits the unit circle as natural boundary
(this was proved by Szeg\H{o} in \cite{Szego}). \endpf

\item[(b)] {\it There exist power series with Hadamard grade $\geq 2$ (i.e., 
non-algebraic sequences with finite Hadamard grade, since all series with
Hadamard grade $\geq 2$ are transcendental).}

\proof There is a classical example of a nonalgebraic Hadamard product of two algebraic
power series (see \cite{Jungen}). Namely the Hadamard square of 
$\sum {2n \choose n} z^n = (1-4z)^{-1/2}$, is equal to $\sum {2n \choose n}^2 z^n = 
\frac{2}{\pi} \int_0^{\pi/2} \frac{\mbox{\rm d} \theta}{(1-16z\sin^2 \theta)^{1/2}}$ 
which is transcendental. \endpf

\item[(c)] {\it There exist $D$-finite power series that have infinite Hadamard grade.
There even exist $D$-finite power series with integer coefficients that have infinite 
Hadamard grade.}

\proof First the exponential series $\sum_{n \geq 0} z^n/n!$ has infinite
Hadamard grade since the denominators of its coefficients are divisible by infinitely 
many primes (actually {\it all} the primes); it is of course $D$-finite. Second the 
Euler series $E(z) := \sum_{n \geq 0} (-1)^n n! z^n$ is clearly $D$-finite 
since $z^2 E'(z) + (1+z) E(z) = 1$, while it has infinite Hadamard grade
(its radius of convergence is $0$). \endpf

{\bf Question.} Is there a $D$-finite power series with integer coefficients and
positive radius of convergence which is algebraic on ${\mathbb F}_p(z)$ for all primes 
$p$ when reduced modulo $p$, but has infinite Hadamard grade?

{\bf Question.} Is it true that algebraic power series can be decomposed into finite
Hadamard products of quadratic power series?

\item[(d)] {\it There exist two power series $A$ and $B$ having infinite Hadamard
grade, such that $A * B$ has finite Hadamard grade.}

\proof Of course if $A$ is even and $B$ odd, then $A * B = 0$. But, from 
Section~\ref{conv}, we have even more than the claim above: take $A(z) := E(z)$ 
the Euler series and $B(z) := e^z$. These power series both have infinite Hadamard 
grade, and their Hadamard product is rational, hence of Hadamard grade $1$. 
Another example is given by the Hadamard square of $\sum \pm z^n$, where
the sequence of $\pm$'s is not ultimately periodic. \endpf

{\bf Questions.} 

\begin{itemize}

\item Is it true that for all $k > 0$ there exists a power series whose grade is $k$?

\item Is it true that if $B$ is a nonzero rational power series, then for every
power series $A$ the Hadamard grade of $A+B$ and of $AB$ are both equal to the 
Hadamard grade of $A$?

\item Is it true that if the Hadamard grade of a power series $A$ is infinite, then for
any nonzero power series $B$ of finite Hadamard grade, the Hadamard grades of $A+B$ and 
$AB$ are both infinite?

\item Let $a$ be an irrational number (real or complex). Is it true that the grade of 
$(1+z)^a$ is infinite?

\item Let $\sum a_n z^n$ be algebraic irrational. Is it true that for all integers
$k \geq 1$ the Hadamard grade of $\sum a_n^k  z^n$ is equal to $k$? Is it the case
if $a_n := {2n \choose n}$?

\item Let $A(z) := \sum a_n z^n$. Let ${\cal H}_k$ denote the Hadamard grade of
the Hadamard product of $k$ factors equal to $A$. If there exists $C > 0$ such 
that ${\cal H}_k < C$ for all $k$, is it true that $A(z)$ must be a rational 
function? What if $A$ satisfies the weaker condition ${\cal H}_k = o(k)$? or 
even the weaker condition that there exists an integer $k \geq 2$ such that 
${\cal H}_k < k$? (Note that conversely if $A$ is rational then ${\cal H}_k = 1$
for all $k \geq 1$.)

\end{itemize}

\end{itemize}

\begin{remark} More examples (where ${\mathbb Q}$ is replaced by ${\mathbb C}$) can
be obtained from results of \cite{BGKL} on non-$D$-finite power series.
\end{remark}

\section{Acknowledgments}
We want to thank warmly J.-P. B\'ezivin, M. Bousquet-M\'elou, and
J.-M. Deshouillers for comments and discussions.

\end{document}